

\def\prop#1;#2\par{\medbreak\noindent{\bf#1\ Proposition.\enskip }{\sl#2}\par
         \ifdim\lastskip<\medskipamount
         \removelastskip\penalty100\medskip\fi}
\def\theorem#1;#2\par{\medbreak\noindent{\bf#1\ Theorem.\enskip }{\sl#2}\par
         \ifdim\lastskip<\medskipamount
         \removelastskip\penalty100\medskip\fi}
\def\lemma#1;#2\par{\medbreak\noindent{\bf#1\ Lemma.\enskip }{\sl#2}\par
         \ifdim\lastskip<\medskipamount
         \removelastskip\penalty100\medskip\fi}
\def\britem#1){\item{#1)}}

\def\proof#1;{\noindent{\bf Proof\setbox0\hbox{#1}%
\ifdim\wd0>.1pt \space \else\fi #1:}\space\ignorespaces}
\def\ni{\noindent}

\def\deq{\mathrel{\lower3pt\hbox{$\buildrel{\rm def}\over=$}}}

\def\nmp{{\rm Nm}_p}

\def\cctil{{\widetilde C}}
\def\dual{^{\vee}}
\def\ctil{\cctil_{\eta}}

\def\ee{{\cal E}}

\def\xieta{{\Xi_{\eta}}}
\def\peta{P_{\eta}}

\def\jtil{{\widetilde J^{2g-2}}}
\def\nm{{\rm Nm}}
\def\pic{{\rm Pic}}
\def\prym{{\rm Prym}}

\def\nn{{\cal N}}
\def\mm{{\cal M}}

\def\ll{{\cal L}}

\def\bull{\item{$\bullet$}}

\def\ll{{\cal L}}

\def\dd{{\cal D}}
\def\space{$$$$$$*$$$$$$}

\def\ses#1#2#3{0\rightarrow {#1}\rightarrow {#2}\rightarrow {#3}\rightarrow 0}

\def\z{{\bf Z}}
\def\ad{{\rm ad}\ }

\def\mm{{\cal M}}
\def\p{{\bf P}}

\def\im{{\rm image \ }}

\def\ker{{\rm ker \ }}

\def\map#1{\ \smash{\mathop{\longrightarrow}\limits^{#1}}\ }

\def\oo{{\cal O}}
\def\pf{{\it Proof.  }}
\def\c{{\bf C}}
\def\th{K^{1\over 2}}
\def\bull{$\bullet$}

\def\c{{\bf C}}
\def\h{{\bf h}}

\def\t{{\bf t}}

\def\z{{\bf Z}}
\def\spin{{\rm Spin}}
\def\pin{{\rm Pin}}
\def\prym{{\rm Prym}}

\def\ll{{\cal L}}
\def\oo{{\cal O}}
\def\dd{{\cal D}}
\def\ee{{\cal E}}

\def\nn{{\cal N}}
\def\mm{{\cal M}}

\def\cotheta{\theta^{\vee}}
\def\zepow#1{(1-\zeta^{#1})}
\def\xipow#1{(1-\xi^{#1})}

\def\ftil{\widetilde F}
\def\gtil{\widetilde g}
\def\pitil{\widetilde \pi}
\def\<{\langle}
\def\>{\rangle}
\def\endo{{\rm End}}
\def\nm{{\rm Nm}}
\def\im{{\rm Im}}

\def\pic{{\rm Pic}}

\def\Po#1;{\Pi_{#1}}
\def\Pd#1;{\Delta_{#1}}
\def\Ps#1;{\Phi_{#1}}
\def\Nd#1;{{\cal Q}_{#1}}
\def\Ns#1;{{\cal N}_{#1}}

\magnification=\magstep1
\font\capit=cmcsc10
\font\addressit=cmcsc8
\font\eightrm=cmr8

\overfullrule=0pt

\centerline{\bf PRYM VARIETIES AND THE MODULI OF SPIN BUNDLES}

\bigskip
\bigskip
\centerline{\bf W.M. Oxbury }
\medskip
\centerline{\it Department of Mathematical Sciences}
\centerline{\it University of Durham}

\medskip

\bigskip
\bigskip

\beginsection Introduction

\medskip

This article is a sequel to the paper [O], where the following construction is
described.
On any smooth projective curve $C$ of genus $g\geq 2$, let $M_d$ (for $d \in \z
/2$) denote the projective moduli space of semistable rank 2 vector bundles
with fixed determinant line bundle of degree $\equiv d$ mod 2. Then there are
natural homomorphisms
 between vector spaces of equal dimension:
$$H^0(M_0,K^{-1})\dual \  \rightarrow \  \bigoplus_{\eta\in
J_2(C)}H^0_+(\peta,3\xieta),\leqno \bf (0.1)
$$
$$H^0(M_1,K^{-1})\dual \  \rightarrow \  \bigoplus_{\eta\in
J_2(C)}H^0_-(\peta,3\xieta);
$$
where $(P_{\eta},\xieta)$ is the principally polarised Prym variety
corresponding to a 2-torsion point $\eta \in J_2(C)$, or is the Jacobian
$(J(C),\Theta)$ in case $\eta =0$; and where $H^0_{\pm}$ denotes even/odd
sections.

The interest of (0.1) is that it is one step in the direction, initiated by
Beauville and others ([B1], [B2], [BNR]), of relating the Verlinde spaces of
$C$---i.e. the vector spaces $H^0(\mm,\ll^{\otimes k})$ where $\mm$ is some
moduli space of vector bun
dles on $C$ and $\ll$ is the ample generator of $\pic(\mm)$ [DN]---to spaces of
classical theta functions associated to the curve. The basic example of this is
the isomorphism
$$
H^0(\mm(n,\oo),\ll)\dual \cong H^0(J(C), n \Theta)\leqno \bf (0.2)
$$
where $\mm(n,\oo)$ is the moduli space of bundles of rank $n$ and trivial
determinant.

The purpose of the present article is to give a generalisation of (0.1) which
ought in principle to describe all Prym theta functions beyond level 3
(at least for {\it odd} level, a restriction that is necessary in the present
work)
in terms of higher moduli spaces, not of vector bundles but of spin bundles.

For any connected reductive algebraic group $G$ there exists a normal
projective moduli variety $\mm(G)$ of semistable principal $G$-bundles on the
curve $C$ ([R2], [KNR]); its connected components are indexed by the
fundamental group $\pi_1(G)$. In the p
articular case $G=SO_m$ there are therefore two components (labelled by the
second Stiefel-Whitney class of the bundles); and each of these is \'etale
covered by a moduli variety of Clifford bundles lifting the $SO_m$-bundles. We
shall denote these two co
vering spaces (defined in \S2 below) by $\mm^{\pm}(\spin_m)$. Then
$\mm^{+}(\spin_m)=\mm(\spin_m) $; whilst for $m=3$ the varieties $ M_d =
\mm^{(-)^d}(\spin_3)$ are precisely the moduli spaces of rank 2 vector bundles
appearing in (0.1).

In the first two sections of this paper we shall review the basic properties of
$\mm(G)$ and of the spaces of spin bundles respectively. Of particular
importance is the construction of the determinant line bundle $\Theta(V)$ over
any $\mm(G)$ associated t
o a finite dimensional representation $G \rightarrow SL(V)$; and in sections 3
and 4 we compute, using the Verlinde formula, the number of sections of the
line bundles $\Theta(\c^m) \rightarrow \mm^{\pm}(\spin_m)$, where $\c^m$
denotes the standard orthog
onal representation.

In \S5 we prove the main result of the paper, which is the construction of
natural homomorphisms, {\it between vector spaces of equal dimension},
$$
H^0(\mm(\spin_{2n+1}),\Theta(\c^{2n+1}))\dual \ \rightarrow \
\bigoplus_{\eta\in J_2(C)}H^0_+(\peta,({2n+1})\xieta),
\leqno \bf (0.3)
$$
$$H^0(\mm^-(\spin_{2n+1}),\Theta(\c^{2n+1}))\dual \  \rightarrow \
\bigoplus_{\eta\in J_2(C)}H^0_-(\peta,({2n+1})\xieta),
$$
for any $n\geq 1$. Of course, the 3-dimensional orthogonal representation of
$\spin_3$ coincides with the adjoint representation, so that $\Theta(\c^3) =
K^{-1}$ and (0.3) reduces to (0.1) in the case $n=1$.

One expects (0.1) and (0.3) to be isomorphisms, at least for the generic curve;
however, we have not been able to prove this. Even more intriguing is the
question of what the corresponding picture is for the even spin groups: not
only do the coincidences
of dimension underlying (0.3) fail, but also our construction of these maps
fails at several points if $m$ is even. But that these groups do indeed have a
tale worth telling is indicated by the tantalising identity for $\spin_8$:
$$
h^0(\mm(\spin_8), \Theta(\c^8)) =
\sum_{\eta \in J_2(C)} h^0(P_{\eta}, 8 \xieta),
$$
(see (5.3)) which one can guess has something to do with triality. On the other
hand, we point out in \S6 that the union $\bigcup_{\eta \in J_2(C)} \peta$ of
all the Prym varities (including the Jacobian) can be interpreted as a moduli
space of $\pin_2$-b
undles over the curve, thus bringing our results into the broader context of
`spin reciprocity', analogous to the `strange duality' of which (0.2) is
supposed to be a special case (see [OW]).

The philosophy which emerges is: in generalising properties of the  rank 2
moduli spaces, and in particular their close relationship with the classical
Schottky-Jung-Prym geometry of the curve, one should pass not only to higher
rank vector bundles, but a
lso to the spin moduli spaces $\mm^{\pm}(\spin_n)$.
For example, the results (0.2) and (0.3) indicate, at the level of theta
functions, a dictionary between abelian varieties and groups:
$$
\eqalign{{\rm Jacobians} \qquad &\leftrightarrow \qquad SL_n \cr
         {\rm Pryms} \qquad &\leftrightarrow \qquad SO_n , \cr}
$$
and this sort of correspondence has appeared also in the work of Kanev [K] in
studying Lax flows for simple Lie algebras and their linearisation on
Prym-Tyurin varieties. It would be fascinating to understand the connection
with this latter point of view.

\medskip
\noindent
{\it Acknowledgements.} In writing this paper the author has benefited from
conversations with M.S. Narasimhan, C. Sorger and E. Viehweg, to whom he
expresses his thanks.

\space

\beginsection \S1 Moduli spaces of $G$-bundles on a curve

\medskip

In this section we shall give a brief account of the moduli spaces of
semistable principal bundles over a curve, following [R1], [DN], [KNR].

We begin with a smooth projective complex curve $C$ of genus $g \geq 2$, and a
complex connected reductive algebraic group $G$; and we consider algebraic
principal $G$-bundles $E\rightarrow C$. Topologically such bundles are
classified by the fundamental
group of $G$. Namely, we can trivialise $E$ over the complement of a point
$x\in C$ and in a small neighbourhood of the point; we then use a loop in $G$
to glue the two trivialisations together on the overlap of the two
neighbourhoods. The homeomorphism c
lass of the resulting bundle depends only on the homotopy class of the loop.

Just as for vector bundles, one has notions of stability, semistability and
S-equivalence for algebraic $G$-bundles, and for stable bundles S-equivalence
is the same as isomorphism. (We shall recall in a moment the definition of
stability, but it will no
t be necessary here to define S-equivalence.)
The basic result of Ramanathan [R2] is then the following.

\proclaim (1.1) Theorem. Given $C,G$ as above and an element $\gamma \in
\pi_1(G)$, there exists a normal irreducible projective variety $\mm(G,\gamma)$
which is a coarse moduli space for families of semistable $G$-bundles on $C$,
modulo S-equivalence.

Moreover, one has
$$
\dim \mm(G,\gamma) = (g-1)\dim G + \dim Z(G),$$
and $\mm(G,\gamma)$ is unirational when $G$ is a simple group [KNR].

The basic construction with principal bundles is the following. If $E$ is a
$G$-bundle, and $\rho:G\rightarrow {\rm Aut}(X)$ any left $G$-space, then we
can form a bundle $E(X) = E\times_{\rho} X$ with fibre $X$. In case $X=G/P$ is
a homogeneous coset spa
ce, a section $\sigma: C \rightarrow E(G/P)$ is called a reduction of the
strucure group of the bundle to the subgroup $P$. When $P\subset G$ is a
maximal parabolic, $E(G/P)\rightarrow C$ can be thought of as
a `generalised Grassmannian bundle'. Then by definition, $E$ is {\it stable}
(resp. {\it semistable}) if and only if
$$
\deg \sigma^* T^{\rm vert}E(G/P) > 0
\quad
({\sl resp.}\ \geq 0)
\qquad
\hbox{\sl for all maximal parabolics $P\subset G$,}
$$
where $T^{\rm vert}$ denotes the vertical tangent bundle.

On the other hand, if $\pi : G' \rightarrow G$ is a group epimorphism then we
can view $X=G$ as a left $G'$-space via $\pi$, and so form a $G$-bundle
$E=F(G)$ from any $G'$-bundle $F$. $F$ is said to be a {\it lift} of $E$. In
particular, if $G'$ is a cen
tral extension of $G$ then there is a bijection between maximal parabolics
$P\subset G$ and maximal parabolics $P'=\pi^{-1} P\subset G'$, and moreover
$F(G'/P') \cong E(G/P)$ if $F$ is any lift of $E$.
Consequently:

\proclaim (1.2) Lemma. If $E$ is a $G$-bundle and $F$ a lift of $E$ to a
central extension of $G$ then $E$ is stable (resp. semistable) if and only if
$F$ is.

Finally, of course, we can take for the $G$-space $X$ a finite-dimensional
representation $\rho:G\rightarrow GL(V)$, to obtain a vector bundle $E(V)$.
In the case when $G = GL_n$ and $V = \c^n$ is the standard representation, the
notions of stability, semistability and S-equivalence are
the same for the principal bundle $E$ as for the vector bundle $E(V)$. Thus
we shall write $\mm(n,d) = \mm(GL_n,d)$, for $d\in \pi_1(GL_n) \cong \z$; this
is then the usual moduli space of semistable vector bundles of rank $n$ and
degree $d$.

Consider now the determinant morphism
$$
\det : \mm(n,d) \rightarrow J^d(C).
$$
One knows from [DN] that, via det, $\mm(n,d)$ has Picard group
$$
\pic \ \mm(n,d) \cong \pic \ J^d(C) \oplus \z\{\Theta_{n,d}\}
$$
where $\Theta_{n,d}$ is an ample line bundle on the fibres
constructed as follows. It will be convenient for our purposes always to assume
that $n$ divides $d$, i.e. that we are dealing with vector bundles of integral
slope.

Consider first an arbitrary family $F \rightarrow C\times S$ of semistable
vector bundles on $C$ with rank $n$, degree $d$ and slope $\mu = d/n \in \z$,
as above; and we construct a line bundle $\Theta (F) \rightarrow S$, functorial
with respect to base c
hange $S' \rightarrow S$, in the following way. Let $\pi : C\times S
\rightarrow S$ be the projection. Then (at least Zariski locally)
there is a homomorphism of locally free sheaves on $S$, $\phi :K^0 \rightarrow
K^1$, having the direct images of $F$ under $\pi$ as kernel and cokernel:
$$
0\rightarrow R^0_{\pi}F
\rightarrow K^0\map{\phi} K^1\rightarrow R^1_{\pi}F\rightarrow 0.
\leqno (1.3)
$$
Moreover, the determinant line bundle
$$
\textstyle
{\rm Det}(F) = (\bigwedge^{\rm top} K^0)\dual \otimes (\bigwedge^{\rm top} K^1)
$$
is well-defined and functorial with respect to base change. If $\mu = g-1$ we
write $\Theta(F) = {\rm Det}(F)$; and this has a canonical section $\det \phi$.
Otherwise $\Theta(F)$ is defined to be a suitable twist of ${\rm Det}(F)$ such
that
$$
\Theta(F\otimes \pi^*L) = \Theta(F)
$$
for any line bundle $L\rightarrow S$; i.e. $\Theta$ respects equivalence of
families.

Now in the case $\mu = g-1$ it is shown in [DN] that the functor $\Theta$ is
represented by a Cartier divisor
$$
\Theta_{n,n(g-1)} = {\rm Closure}\{{\rm stable}\ V| H^0(V) = H^1(V) \not= 0 \}
\subset \mm(n,n(g-1)).
$$
In other words $\Theta(F) = f^* \Theta_{n,n(g-1)}$ where $f: S \rightarrow
\mm(n,n(g-1))$ is given by the coarse moduli property.

For the general case one chooses a line bundle $L\in \pic(C)$, with integral
slope, and degree chosen so that we get a morphism
$$
\mm(n,d) \map{\otimes L} \mm(n,n(g-1)).
$$
Now set
$$
\Theta_L = (\otimes L)^*\Theta_{n,n(g-1)}.
$$
The dependence of $\Theta_L$ on $L$ is then given by (1.5) below:

\proclaim (1.4) Lemma [DN].
View $J^0(C)$ as a subgroup of $\pic(J^d(C))$ by $L \mapsto \Phi^{-1}\otimes
T^*_L\Phi$, where $\Phi$ is any line bundle representing the principal
polarisation on $J^d(C)$. Then for any family $F\rightarrow C\times S$ as
above, and any $L\in J^0(C)$, we
have
$$
\Theta(F\otimes pr_C^*L) = {\det}^*(L) \otimes \Theta(F)
$$
where $\det : S\rightarrow \mm(n,d) \rightarrow J^d(C)$.

It follows easily from this that when $L,L'$ have the same degree,
$\Theta_L$ and $\Theta_{L'}$ are related by:
$$
\Theta_{L'} = {\det} ^*(L'L^{-1}) \otimes \Theta_L \in \pic\ \mm(n,d).
\leqno(1.5)
$$
We now set $\Theta_{n,d} = \Theta _L$: if $d\not= n(g-1) $ this depends on $L$,
{\it but by (1.5) its restriction to the fibres of $\det :\mm(n,d) \rightarrow
J^d(C)$
is independent of $L$}.

We need to remark that in fact a more general form of (1.5) is given in [DN].
Namely, $L$ (and $L'$) can be a vector bundle of higher rank  chosen so that
tensoring $\mm(n,d)$ by $L$ gives bundles of slope $g-1$; so that then
$\Theta_L$ can be defined in
the same way as above. In this more general situation (1.5) becomes:
$$
\Theta_{L'} = {\det} ^*(\det L' \otimes \det L\dual) \otimes \Theta_L \in \pic\
\mm(n,d).
\leqno(1.5')
$$

\medskip

Now suppose that we are given a family $E\rightarrow C\times S$ of semistable
$G$-bundles, and a representation $\rho:G\rightarrow SL(V)$, where $\dim V =n$.
Then we can form the family of vector bundles $E(V) \rightarrow C\times S$;
that these are semistable follows from [R2] proposition 2.17.
So we have a theta line bundle $\Theta(E(V))\rightarrow S$, and since $E(V)$
has trivial determinant on the fibres of $\pi:C\times S\rightarrow S$ we deduce
from lemma (1.4) the following corollary, which will be needed later:

\proclaim (1.6) Corollary. For $E\rightarrow C\times S$ and  $\rho:G\rightarrow
SL(V)$ as above, and for any $L\in J^0(C)$ one has
$$
\Theta(E(V)\otimes pr_C^*L) = \Theta(E(V)).
$$

Globally $\rho$ induces a morphism
$$
\rho_* : \mm(G,\gamma) \rightarrow \mm(SL_n) \hookrightarrow \mm(n,0)
$$
and the functor $E\mapsto \Theta(E(V))$ is represented by the (well-defined)
line bundle
$$
\Theta(V) := (\rho_*)^* \Theta_{n,0} \in \pic \ \mm(G,\gamma).
$$

Note that if we let $j: \mm(SL_n) \hookrightarrow \mm(n,0)$ denote the
inclusion which identifies $\mm(SL_n)$
with the moduli space of vector bundles of rank $n$ and trivial determinant,
via the standard representation $\c^n$, then by construction
$$
\Theta(\c^n) = j^*\Theta_{n,0}.
$$

\space

\beginsection \S2 Clifford bundles

\medskip

Let us consider again the fibration
$$
\eqalign{
\det :\ &\mm(n,d) \rightarrow J^d(C);\cr
      }
$$
induced, that is, by
the determinant homomorphism $GL_n \rightarrow \c^*$.
Writing $M_d$ for the isomorphism class of the fibre, $M_d$ depends on $d$ only
modulo $n$, and so the structure of $\mm(GL_n)$ reduces to that of a set of $n$
varieties, of which that corresponding to $d\equiv 0$ mod $n$
is isomorphic to $\mm(SL_n)$.

In this section we shall describe an alternative generalisation of this
situation for $n=2$, obtained by replacing $GL_2$ not by $GL_n$, but
by the special Clifford group of a nondegenerate quadratic form. First we need
to recall briefly some basic Clifford theory. (See [Bo].)

Let $Q$ be a nondegenerate quadratic form on a complex vector space $V$ of
finite dimension $m$;
let $A=A(Q)$ be its Clifford algebra and $A^+$ the even Clifford algebra.
Recall that these
can be expressed as matrix algebras as follows.

If $m=2n$ is even then for any $n$-dimensional isotropic subspace $U\subset V$
one has
$$
\textstyle
A \cong \endo \bigwedge U;
\qquad
A^+ \cong \endo \bigl( \bigwedge^{\rm even} U \bigr) \oplus \endo \bigl(
\bigwedge^{\rm odd} U \bigr).
$$

If, on the other hand, $m=2n+1$ is odd then for any direct sum decomposition
$V=U \oplus U'\oplus \c$ where $U,U'$ are $n$-dimensional isotropic subspaces
one has
$$
\textstyle
A \cong \endo \bigwedge U \oplus \endo \bigwedge U';
\qquad
A^+ \cong \endo \bigwedge U .
$$
The `principal involution' of $A$ is $\alpha : x\mapsto -x$ for $x\in V$, i.e.
is $\pm 1$ on $A^{\pm}$ respectively. The `principal anti-involution' $\beta$
is the identity on $V$ and reverses the direction of multiplication:
$\beta(x_1\ldots x_r) = x_r\ldots x_1$. Then {\it conjugation} is defined by
$$
x^* = \beta \circ \alpha (x),
\qquad x\in A.
$$
Finally one can define the {\it Clifford group}
$$
 C (Q) = \{s\in  A^{*} | sVs^* \subset V\},
$$
where $A^{*}\subset A$ denotes the group of units; and the {\it special
Clifford group}
$$
SC(Q) = C(Q) \cap A^+.
$$
For $s\in C(Q)$ the transformation $\pi_s: x\mapsto sxs^*$ of $V$ is
orthogonal---this is because $C(Q)$ is generated by $x\in V\cap C(Q)$, for
which $\pi_x$ is just reflection in the hyperplane $x^{\perp}$. Thus one has a
group homomorphism $\pi:C(Q) \rightarrow O(Q)$, which
has the following properties.

\proclaim (2.1) Proposition.
\item{(i)} $\ker \pi = C(Q) \cap Z(A)$;
\item{(ii)} if $m$ is odd then $\pi(C(Q)) = \pi(SC(Q)) = SO(Q)$;
\item{(iii)} if $m$ is even then $\pi(C(Q)) = O(Q)$ and $\pi(SC(Q)) = SO(Q)$.

\proclaim (2.2) Corollary. $SC(Q)$ is a connected reductive algebraic group.

The {\it spinor norm} is the group homomorphism
$$
\eqalign{
\nm \ :\ C(Q) &\rightarrow \c^{*} \cr
         s &\mapsto \beta(s)s.\cr}
$$
Equivalently $\nm (x_1\ldots x_r) =  Q(x_1)\cdots Q(x_r)$ for $x_1,\ldots ,x_r
\in V$. Then $\spin (Q)$ is by definition $SC(Q) \cap \ker \nm$.

 From now on we shall write $C_m, SC_m, \spin_m$ instead of
$C(Q),SC(Q), \spin(Q)$ when $Q$ is the standard quadratic form on $\c^m$.
Then
(using (2.1)) one has the following commutative diagram of short exact
sequences:
$$
\matrix{
1&\rightarrow&\c^{*}&\rightarrow&SC_m&\map{\pi}&SO_m&\rightarrow&1\cr
&&&&&&&&\cr
&&\uparrow&&\uparrow&&\Vert&&\cr
&&&&&&&&\cr
1&\rightarrow&\z/2&\rightarrow&\spin_m&\rightarrow&SO_m&\rightarrow&1.\cr}
\leqno(2.3)$$

\proclaim (2.4) Proposition.
\item{(i)} $SC_m$ has centre $Z(SC_m) =
\cases{\c^{*} & if $m$ is odd \cr
       \c^{*} \times \z/2 & if $m$ is even;\cr}$
\item{(ii)} $SC_m$ has fundamental group $\pi_1(SC_m) = \z$; and this maps
isomorphically to $\pi_1(\c^*) = \z$ under the spinor norm.

\pf (i) If $m$ is odd the centre of $SC_m$ must be contained in---and hence
equal to---the kernel $\c^*$ of the surjection onto $SO_m$, since the latter
has trivial centre.

If, on the other hand,
 $m$ is even, then by the same token $Z(SC_m)$ is contained in
$\pi^{-1}(Z(SO_m))$ where $\pi$ denotes the surjection to $SO_m$. In this case
$SO_m$ has centre $\{\pm 1\}$. As before everything in $\pi^{-1}(1) = \c^*$
is central; while if $\{ e_1, \ldots , e_m \} \subset \c^m$ is any orthonormal
basis then the product $e_1 \ldots  e_m \in SC_m$ spans $\pi^{-1}(-1) \cong
\c^*$. Since $m$ is even this product anticommutes with each $e_i$, and
therefore {\it commutes} wit
h all elements of $A^+$. So $\{\pm 1\} \cong \c^{*} \times \z/2$ is contained
in and therefore equal to the centre.

(ii) From the exact homotopy sequence of the fibration in the upper sequence of
(2.3), and the vanishing of $\pi_2(\c^*)$, we have a non-split extension
$$
\ses{\z}{\pi_1(SC_m)}{\z /2}.
$$
Since the fundamental group of a Lie group is abelian it follows that the only
possibility is $\pi_1(SC_m) = \z$. The last part now follows from the fact that
$\spin_m$ is simply-connected.

\medskip

\noindent
{\bf Example.}
When $m=3$, $SC_3$ is equal to the group of units in $A^+$, i.e. $SC_3 \cong
GL_2$. The homomorphism onto $SO_3$ in (2.3) is then precisely the action of
$GL_2$ on $S^2 \subset {\bf R}^2 \subset \c^3$ by M\"obius transformations, via
stereographic project
ion.
The spinor norm is in this case the determinant homomorphism, so that $\spin_3
= SL_2$.

\medskip

 From (2.4)(ii) we see that there is for each $d\in \z = \pi_1(SC_m)$ a
morphism
induced by the spinor norm, and which we shall denote in the same way:
$$
\nm : \mm (SC_m, d) \rightarrow J^d(C).
$$
Moreover, when $m=3$ this is nothing but the determinant morphism for rank 2
vector bundles. And just as for rank 2 vector bundles, one has:

\proclaim (2.5) Proposition. For $d\in \z$ and $L\in J^d(C)$, the isomorphism
class of the subscheme $\nm ^{-1}(L) \subset \mm(SC_m,d)$ depends only on $d$
mod 2.

This is a trivial result proved in the same way as for rank 2 vector bundles,
once one observes that multiplication of $SC_m$ by its centre (2.4)(i) induces
a natural generalisation of the
tensor product operation of Clifford bundles by line bundles if $m$ is odd, and
if $m$ is even by {\it pairs} $(N,\eta)$ where $N$ is a line bundle and $\eta
\in H^1(C,\z/2) = J_2(C)$. (See also [R].)
We shall write, respectively, $N\otimes E$ and $(N,\eta)\otimes E$ for this
product.
It follows from the definition of the spinor norm that
$$
\nm(N\otimes E) = N^2 \otimes \nm (E),
\quad {\rm resp.}
\quad \nm((N,\eta)\otimes E) = N^2 \otimes \nm (E).
$$

So to prove (2.5): suppose $d=\deg L \equiv d'=\deg L'$ mod 2, and write $L' =
N^{2k}\otimes L$ for some $N\in \pic^1(C)$, $k\in \z$. Then the map
$\mm(SC_m,d)
\rightarrow \mm(SC_m,d')$ given by $E\mapsto N\otimes E$ (resp. $(N,\oo)\otimes
E$) restricts to an isomorphism $\nm^{-1}(L) \cong \nm^{-1}(L')$.

We shall therefore introduce the notation:
$$
\eqalign{
\mm^+(\spin_m) &= \mm(\spin_m) = \nm^{-1}(\oo_C);\cr
\mm^-(\spin_m) &= \nm^{-1}(\oo_C(p)),\cr}
$$
where $p\in C$ is any point of the curve. The second space generalises the
moduli space of stable rank 2 vector bundles with fixed determinant of odd
degree, in the case $m=3$.

We next want to relate the two varieties $\mm^{\pm}(\spin_m)$ to the
two-component moduli variety $\mm(SO_m)$ of semistable $SO_m$-bundles. Given
any $SC_n$-bundle $E \rightarrow C$ we can construct an $SO_m$-bundle
$E(SO_m)$ as in \S1, using the surjection $\pi$ in (2.3). Then the following is
a special case of lemma (1.2).

\proclaim (2.6) Lemma. An $SC_m$-bundle $E$ is stable (resp. semistable) if and
only if $F=E(SO_m)$ is stable (resp. semistable).

\medskip

For an algebraic group $G$, let $\oo G$ denote the sheaf of holomorphic
functions on $C$ with values in $G$. Then isomorphism classes of $G$-bundles on
$C$ are parametrised by the Cech cohomology group
$H^1(C,\oo G)$. Diagram (2.3) together with the spinor norm induces the
following exact commutative diagram of such cohomology groups (see also [R]):
$$
\matrix{{\bf (2.7)}
&&\pic(C)&=&\pic (C)&&&&\cr
&&&&&&&&\cr
&&\uparrow sq&&\ \ \uparrow \nm&&&&\cr
&&&&&&&&\cr
0&\rightarrow&\pic (C)&\rightarrow&H^1(C,\oo SC_m)&\rightarrow&H^1(C,\oo
SO_m)&\rightarrow&H^2(C,\oo^*) =0\cr
&&&&&&&&\cr
&&\uparrow&&\uparrow&&\parallel&&\cr
&&&&&&&&\cr
0&\rightarrow&J_2(C)&\rightarrow&H^1(C,\oo \spin_m)&\rightarrow&H^1(C,\oo
SO_m)&\map{w_2}&H^2(C,\z/2) = \z/2.\cr
}
$$
Here $sq$ is the squaring map, and the map $H^1(C,\oo SC_m)\rightarrow
H^1(C,\oo SO_m)$ takes the isomorphism class of $E$ to that of $E(SO_m)$.
The connecting homomorphism
$w_2$ is the second Stiefel-Whitney class, and it can be seen from (2.7) that
$$
w_2(E(SO_m)) \equiv \deg \nm (E) \quad {\rm mod} \ 2.
\leqno(2.8)
$$
Note also that whereas the vanishing of $w_2$ is the necessary and sufficient
condition for an $SO_m$-bundle to lift to a $\spin_m$-bundle, {\it every}
$SO_m$-bundle lifts to a $SC_m$-bundle, the space of such lifts being bijective
to $\pic(C)$. Moreover,
 one has:

\proclaim (2.9) Proposition. Given an $SO_m$-bundle $F$ and a line bundle $L\in
\pic(C)$, there is a (non-canonical) bijection between the set of
$SC_m$-bundles $E$ lifting $F$ with $\nm(E) = L$, and $J_2(C)$.

The group $SO_m$ has fundamental group $\z/2$; it will be convenient to write
$\mm(SO_m,0) = \mm^+(SO_m)$ and $\mm(SO_m,1) = \mm^-(SO_m)$. Then it follows
from (2.6), (2.8) and (2.9) that there are natural Galois covers
$$
\mm^{\pm}(\spin_m) \map{J_2(C)} \mm^{\pm}(SO_m).
\leqno (2.10)
$$

Finally, we shall need the following fact relating Clifford bundles with
theta-characteristics.

\proclaim (2.11) Proposition. For any $SC_m$-bundle $E\rightarrow C$ and
theta-characteristic $\th \in J^{g-1}(C)$ one has
$$
h^0(C,\th \otimes E(\c^m)) \equiv m h^0(C,\th) + \deg \nm(E)
\quad {\rm mod}\ 2.
$$

{\it Proof.} For any orthogonal bundle $F$ and theta-characteristic $\th$, one
has by [S], th\'eor\`eme 2, the congruence:
$$
h^0(C,\th \otimes F) \equiv (m+1) h^0(C,\th) + h^0(C,\th \otimes w_1(F))+
w_2(F)
\quad {\rm mod}\ 2,
$$
where $w_1$ and $w_2$ are the Stiefel-Whitney classes. But $w_1(F)$ can be
identified with $\det F \in J_2(C) \cong H^1(C,\z/2)$, which in the present
case $F=E(\c^m)$ vanishes since $E$ is a {\it special} Clifford bundle. On the
other hand, $w_2(F)\equiv
 \deg \nm(E)$ mod 2 by (2.8). So we get the statement in the proposition.

\space

\beginsection \S3 The Verlinde formulae

\medskip

In this section we shall write down, for the unitary and spin groups, the
Verlinde formula which calculates the dimension of the vector spaces
$H^0(\mm(G),\Theta(V))$.
For the computations of these formulae we refer the reader to [B3], [OW].
In fact, what one writes down is a natural number $N_l(G)$ (or to all
appearances a somewhat unnatural number!)
depending on the group, on the genus $g$, and on an integer $l$ called the
`level'. Then in the next section we shall recall how, to any representation
$V$ of $G$, one associates a level $l=m_V$---its `height'---such that
$\dim H^0(\mm(G),\Theta(V)) = N_{m_V}(G)$.

\space

\noindent $\bf SL_{n+1}$

\medskip

The Verlinde number in this case is:
$$
\eqalign{
N_l(SL_{n+1}) =& \bigl((n+1)(l+n+1)^n\bigr)^{g-1}\times \cr
              & \sum_{\scriptstyle  t_1,\ldots ,t_n \geq 1
                     \atop \scriptstyle t_1 +\cdots +t_n \leq l+n}
     \Bigl( \prod_{1\leq i<j\leq n+1}
           \bigl(2 \sin {\pi \over l+n+1}(t_i +\cdots +t_{j-1}) \bigr)
       \Bigr)^{-2g+2}.\cr}
\leqno(3.1)
$$

Since
we are especially concerned with the spin groups we shall record at this point
the result of computing (3.1) at levels 1 and 2 for $\spin_6 = SL_4$:
\item{\bull} $N_1(SL_4) = 2^{2g}$;
\item{\bull} $N_2(SL_4) = 2^{3g-1}3^{g-1} + 2^{3g-1} + 2^g 3^{g-1}$.

\noindent Also for reference below let us list here the well-known numbers for
$SL_2$:
\item{\bull} $N_1(SL_2) = 2^{g}$;
\item{\bull} $N_2(SL_2) = 2^{g-1}(2^g +1)$.

\bigskip

\noindent $\bf Spin_{2n},\ n \geq 4$

\medskip

Here the Verlinde number is:
$$
N_l(\spin_{2n})= \bigl(4(l+2n-2)^n\bigr)^{g-1} \times
                  \sum_{\scriptstyle
                      t_1,\ldots,t_n \geq 1 \atop \scriptstyle
                         t_1 +2t_2 +\cdots 2t_{n-2} +t_{n-1}+t_n \leq l+2n-3}
                     \Bigl( \prod_{i=1}^{n-1} \times \prod_{i<j}
                      \Bigr)^{-2g+2},
\leqno(3.2)
$$
where
$$
\prod_{i=1}^{n-1} =  \prod_{i=1}^{n-1} 4 \sin {\pi\over l+2n-2}(t_i +\cdots
+t_{n-1})
                                  \sin {\pi\over l+2n-2}(t_i +\cdots
+t_{n-2}+t_n),
$$
and
$$
\eqalign{
\prod_{i<j} =& \prod_{1\leq i<j\leq n-1}
4 \sin {\pi\over l+2n-2}(t_i+\cdots +t_{j-1}) \times \cr
 & \sin {\pi\over l+2n-2}(t_i+\cdots +t_{j-1}+2t_j+\cdots +2t_{n-2}
+t_{n-1}+t_n).\cr}
$$
Again we note the first two cases of this formula:
\medskip
\item{\bull} $N_1(\spin_{2n}) = 2^{2g}$;
\item{\bull} $\eqalign{
              N_2(\spin_{2n}) 
                             &= (2n)^g +(2^{2g} -1)(2n)^{g-1} + 2^{g-1}(2^{2g}
                                 -n^g).\cr}$
\medskip

The derivations of these and similar formulae for $\spin_{2n+1}$ below are all
entirely similar (and rather tedious) and we shall content ourselves with
giving it only for the case relevant to the present paper,
$N_2(\spin_{2n+1})$ below.

Note, incidentally, that both of the above formulae are consistent with
$\spin_4 = SL_2\times SL_2$ and $\spin_6 = SL_4$: for $n=2,3$ respectively they
coincide with $N_l(SL_2)^2$ and $N_l(SL_4)$.

\bigskip

\noindent $\bf Spin_{2n+1}, n \geq 2$

\medskip

Here we have:
$$
N_l(\spin_{2n+1})= \bigl(4(l+2n-1)^n\bigr)^{g-1} \times
                  \sum_{\scriptstyle
                           t_1,\ldots,t_n \geq 1 \atop \scriptstyle
                         t_1 +2t_2 +\cdots 2t_{n-1} +t_n \leq l+2n-2}
                     \Bigl( \prod_{i=1}^n \times \prod_{i<j}
                      \Bigr)^{-2g+2},
\leqno(3.3)
$$
where
$$
\prod_{i=1}^n =  \prod_{i=1}^n 2 \sin {\pi\over l+2n-1}(t_i +\cdots +t_{n-1}
+{t_n\over 2}),
$$
and
$$
\eqalign{
\prod_{i<j} =& \prod_{1\leq i<j\leq n}
4 \sin {\pi\over l+2n-1}(t_i+\cdots +t_{j-1})\times  \cr
 & \sin {\pi\over l+2n-1}(t_i+\cdots +t_{j-1}+2t_j+\cdots +2t_{n-1} +t_n).\cr}
$$

This time we obtain at levels 1 and 2:
\medskip
\item{\bull} $N_1(\spin_{2n+1}) = 2^{g-1}(2^g +1)$;
\item{\bull} $N_2(\spin_{2n+1}) = 2^{2g-1}(2n+1)^{g-1} + n(2n+1)^{g-1}
+2^{2g-1}$.
\medskip

We shall prove the second of these formulae. As already mentioned, this is
entirely representative of the proofs of all similar special cases of the
Verlinde formulae given in this section.

\bigskip

\noindent
{\bf Calculation of $N_2(\spin_{2n+1})$.}

\medskip

Let $\zeta_k = e^{2\pi i /k} $; in particular we shall be concerned with $\xi =
\zeta_{4n+2}$ and $\zeta =\xi^2 = \zeta_{2n+1}$. Using the fact that for any
integer $a$,
$$
4\sin^2 {\pi \over k}a = (1- \zeta_k^a)(1-\zeta_k^{-a}),
$$
the Verlinde formula can be written:
$$
N_2(\spin_{2n+1}) = \bigl( 4(2n+1)^n\bigr)^{g-1} \sum_{\scriptstyle
                    t_1,\ldots,t_n \geq 1 \atop \scriptstyle
                         t_1 +2t_2 +\cdots 2t_{n-1} +t_n \leq 2n}(A_{\t}
B_{\t})^{-g+1}
$$
where for $\t = (t_1,\ldots,t_n)$,
$$
\eqalign{
A_{\t} &=  \prod_{i=1}^n (1- \xi^{2t_i+\cdots
+2t_{n-1}+t_n})(1-\xi^{4n+2-2t_i-\cdots -2t_{n-1}-t_n}),\cr
B_{\t} &= \prod_{i<j}(1-\zeta^{t_i+\cdots +t_{j-1}})
         (1-\zeta^{2n+1-t_i-\cdots -t_{j-1}})\times\cr
       & \ \ \ \ \ \ \ (1-\zeta^{t_i+\cdots +t_{j-1}+2t_j+\cdots
+2t_{n-1}+t_n})
       (1-\zeta^{2n+1-t_i-\cdots -t_{j-1}-2t_j-\cdots -2t_{n-1}-t_n}).\cr}
$$
In this sum $\t$ can take only a small number of values---namely, since
$$
2n-2 \leq t_1 +2t_2 +\cdots 2t_{n-1} +t_n \leq 2n,
$$
$\t$ must be one of
$$
\eqalign{
& (1,\ldots,1),\cr
& (2,1,\ldots,1,2),(1,\ldots,1,3),(3,1,\ldots,1),\cr
& (2,1,\ldots,1),(1,2,1,\ldots,1),\ldots ,(1,\ldots,1,2).\cr}
$$
So we just carry out the (somewhat laborious) task of computing $A_{\t}$ and
$B_{\t}$ for each of these vectors. In doing so we shall make repeated use of
the easy identities
$$
\prod_{i=1}^{4n+1} (1-\xi^i) = 4n+2,
\qquad
\prod_{i=1}^{2n} (1-\zeta^i) = 2n+1.
$$

\medskip
\noindent
{\it The case $\t = (1,\ldots,1)$.}

First,
$$
\eqalign{
A_{(1,\ldots,1)} &= \prod_{i=1}^n (1-\xi^{2n-2i+1})(1-\xi^{2n+2i+1})\cr
     &= (1-\xi)(1-\xi^3)\cdots
(1-\xi^{2n-1})(1-\xi^{2n+3})\cdots(1-\xi^{4n+1})\cr
   &=\prod_{i=1}^{4n+1} \xipow{i} /\Bigl( \xipow{2n+1} \prod_{i=1}^{2n}
\xipow{2i}\Bigr) \cr
    &= 1.\cr}
$$
Next, we can write
$$
\eqalign{
B_{(1,\ldots,1)} &= \prod_{i<j}
\zepow{j-i}\zepow{2n+1-j+i}\zepow{2n+1-j-i}\zepow{i+j}\cr
        &= \prod_{k=1}^{2n} \zepow{k}^{m_k}\cr}
$$
where we have to determine the multiplicities $m_k$. Since the factors all
occur in conjugate pairs we have $m_k = m_{2n+1-k}$, so it is sufficient to
assume $1\leq k \leq n$.
Note that the four exponents occuring in the product satisfy the inequalities
$$
\eqalign{
1 \leq & j-i \leq n-1, \cr
n+2 \leq & 2n+1-j+i \leq 2n, \cr
2 \leq & 2n+1-j-i \leq 2n-2, \cr
3 \leq & i+j \leq 2n-1.\cr}
$$
Thus, for $1\leq k \leq n$,
$$
\eqalign{
m_k =& \  \sharp\ \{\hbox{pairs $i<j$ such that $j-i=k$}\}\cr
     & \ +\ \sharp\ \{\hbox{pairs $i<j$ such that $j+i=k$}\}\cr
     & \ +\ \sharp\ \{\hbox{pairs $i<j$ such that $j+i=2n+1-k$}\}\cr}
$$
and the last term is equal to the number of pairs such that $i+j=k-1$. So
$$
\eqalign{
m_k &= n-k + [{k-1\over 2}] + [{k-2\over 2}] \cr
    &= n-k +k-1\cr
    &=n-1,\cr}
$$
for all $k$. Hence
$$
B_{(1,\ldots,1)} = \Bigl( \prod_{i=1}^{2n} \zepow{i} \Bigr)^{n-1} =
(2n+1)^{n-1}.
$$

\medskip
\noindent
{\it The case $\t = (1,\ldots,1,2)$.}

We shall deal with this and the remaining cases by comparing them with the
previous one. So, for example, looking at the definition of $A_{\t}$ we see
that to obtain $A_{(1,\ldots,1,2)}$ we apply the following procedure to
$A_{(1,\ldots,1)}$:
$$
A_{(1,\ldots,1)} = \underbrace{
(1-\xi)(1-\xi^3)\cdots (1-\xi^{2n-1})}_{\scriptstyle\rm add\  1\ to\ each\
exponent}
\underbrace{
(1-\xi^{2n+3})\cdots(1-\xi^{4n+1})}_{\scriptstyle\rm
subtract\ 1\ from\ each\ exponent}
$$
This yields
$$
A_{(1,\ldots,1,2)} = \prod_{i=1}^{2n} \xipow{2i} = \prod_{i=1}^{2n} \zepow{i} =
2n+1.
$$
Similarly, to obtain $B_{(1,\ldots,1,2)}$ we operate on the exponents of
$B_{(1,\ldots,1)}$ by:
$$
B_{(1,\ldots,1)} = \prod_{i<j}
\underbrace{
\zepow{j-i}\zepow{2n+1-j+i}}_{\scriptstyle\rm unchanged}
\underbrace{
\zepow{2n+1-j-i}}_{\scriptstyle\rm add\ 1}
\underbrace{\zepow{i+j}}_{\scriptstyle\rm subtract\ 1}
$$
We claim that when this is done the multiplicity of each factor $\zepow{k}$ in
the resulting product is the same as in $B_{(1,\ldots,1)}$, and hence that
$$
B_{(1,\ldots,1,2)} = B_{(1,\ldots,1)} = (2n+1)^{n-1}.
$$

To prove this, first fix $i\geq 1$. Then one sees that the effect of the above
operation on all the factors coming from $j=i+1,\ldots,n$ is to multiply the
original product by a coefficient
$$
\mu_i = {\zepow{2n-2i+1}\zepow{2i} \over \zepow{n+1-i}\zepow{i+n}}.
$$
Consequently $B_{(1,\ldots,1,2)} = \mu \times B_{(1,\ldots,1)}$ where
$$
\eqalign{
\mu &= \prod_{i=1}^{n-1} \mu_i \cr
    &= {\zepow{2n-1}\cdots \hbox{\it odd powers} \cdots\zepow{3}\zepow{2}
       \cdots \hbox{\it even powers} \cdots \zepow{2n-2} \over
        \zepow{n}\cdots \zepow{2}\zepow{n+1}\cdots \zepow{2n-1}} \cr
    &= 1.\cr}
$$

\medskip
\noindent
{\it The cases $\t = (1,\ldots,1,2,1,\ldots,1)$, $t_k =2$ for $1\leq k\leq
n-1$.}

If we apply the same reasoning as in the previous case then
$A_{(1,\ldots,2,\ldots,1)}$ is obtained from $A_{(1,\ldots,1)}$ by the scheme:
$$
A_{(1,\ldots,1)} =
\underbrace{\xipow{1}\cdots \xipow{2n-2k-1}}_{\scriptstyle\rm unchanged}
\underbrace{\xipow{2n-2k+1}\cdots \xipow{2n-1}}_{\scriptstyle\rm add\ 2}.
$$$$
\times
\underbrace{\xipow{2n+3}\cdots \xipow{2n+2k+1}}_{\scriptstyle\rm subtract\ 2}
\underbrace{\xipow{2n+2k+3}\cdots \xipow{4n+1}}_{\scriptstyle\rm unchanged}
$$
This has the effect of multiplying $ A_{(1,\ldots,1)} $ by the expression
(involving only odd powers of $\xi$)
$$
{\xipow{2n-2k+3}\cdots \xipow{2n+1} \over \xipow{2n-2k+1}\cdots \xipow{2n-1}}
{\xipow{2n+1}\cdots \xipow{2n+2k-1} \over \xipow{2n+3}\cdots \xipow{2n+2k+1}}
$$$$
= {\xipow{2n+1} \xipow{2n+1} \over \xipow{2n-2k+1}\xipow{2n+2k+1}};
$$
i.e. we have
$$
A_{(1,\ldots,2,\ldots,1)} = {4\over \xipow{2n-2k+1}\xipow{2n+2k+1}}
A_{(1,\ldots,1)}.
$$

Next we consider $B_{(1,\ldots,2,\ldots,1)}$. First write
$$
B_{(1,\ldots,1)} = \prod_{k<i<j} \prod_{i\leq k<j} \prod_{i<j\leq k}
\zepow{j-i}\zepow{2n+1-j+i}\zepow{2n+1-j-i}\zepow{i+j}.
$$
To obtain $B_{(1,\ldots,2,\ldots,1)}$ the procedure is:
$$
\underbrace{\zepow{j-i}}_{\scriptstyle\rm add\ 1\ \atop \scriptstyle in\ i\leq
k<j\ only}
\underbrace{\zepow{2n+1-j+i}}_{\scriptstyle\rm subtract\ 1\ \atop \scriptstyle
in\ i\leq k<j\ only}
\underbrace{\zepow{2n+1-j-i}}_{{\scriptstyle\rm add\ 1\ \atop \scriptstyle in\
i\leq k<j;} \atop {\scriptstyle\rm add\ 2\ in\atop \scriptstyle i<j\leq k }}
\underbrace{\zepow{i+j}}_{{\scriptstyle\rm subtract\ 1\ \atop \scriptstyle in\
i\leq k<j;} \atop {\scriptstyle\rm subtract\ 2\ in\atop \scriptstyle i<j\leq k
}}
$$
Now fix $i\leq k$. It can be seen that when we pass to
$B_{(1,\ldots,2,\ldots,1)}$ the factors in the above product coming from
$j=i+1, \ldots ,n$ yield a multiplier
$$
\mu_i = {\zepow{2n-2i+2}\zepow{2n-2i+1}\zepow{2i-1}\zepow{2i}\over
\zepow{k+1-i}\zepow{2n-k+i}\zepow{2n+2-k-i}\zepow{i+k-1}},
$$
and hence that
$$
B_{(1,\ldots,2,\ldots,1)} = \mu \times B_{(1,\ldots,1)}
$$
where
$$
\mu = \prod_{i=1}^k \mu_i =
{\zepow{2k}\zepow{2n+1-2k}\over\zepow{k}\zepow{2n+1-k}}.
$$

We conclude (since $\zeta = \xi^2$; and using also $\xi^{2n+1} = -1$) that
$$
\eqalign{
{A_{(1,\ldots,2,\ldots,1)}B_{(1,\ldots,2,\ldots,1)}\over
A_{(1,\ldots,1)}B_{(1,\ldots,1)}} &= {4\xipow{4k}\xipow{4n+2-4k}\over
\xipow{2n-2k+1}\xipow{2n+2k+1}
\xipow{2k}\xipow{4n+2-2k}}\cr
&= {4 |1-\xi^{4k} |^2 \over |1+\xi^{2k}|^2 |1-\xi^{2k}|^2} \cr
&= 4.\cr}
$$
Hence
$$
A_{(1,\ldots,2,\ldots,1)}B_{(1,\ldots,2,\ldots,1)} = 4(2n+1)^{n-1}.
$$

\medskip

The remaining three cases involve no new ideas and the results are as follows.
The integer values of $A_{\t}B_{\t}$ are:

\medskip
\settabs 4 \columns
\+ &$(2n+1)^{n-1}$ & for $(1,\ldots,1),(3,1,\ldots,1)$;  &\cr
\+ &$(2n+1)^n$ & for $(1,\ldots,1,2),(2,1,\ldots,1,2)$; &\cr
\+ &$4(2n+1)^{n-1}$ & for all other $\t$. &\cr
\medskip
So we arrive at
$$
\eqalign{
N_2(\spin_{2n+1}) &= \bigl( 4(2n+1)^n \bigr)
     \Bigl( {2\over (2n+1)^{(n-1)(g-1)}} + {2\over (2n+1)^{n(g-1)}}\cr
&\ \ \ \ \ \ \  + {n\over (4(2n+1)^{n-1})^{g-1}} \Bigr) \cr
&=  2^{2g-1}(2n+1)^{g-1} + n(2n+1)^{g-1} +2^{2g-1}.\cr}
$$

Finally, it will be convenient, following [OW], to split the Verlinde number
(3.3) for the odd spin groups (including $\spin_3$) into the sum taken over
highest weights of tensor representations---i.e. those which descend to
$SO_{2n+1}$, which is equivale
nt to the condition that $t_n$ is even---and the sum taken over highest weights
of spinor representations, i.e. those for which $t_n$ is odd:
$$
N_l(\spin_{2n+1}) = N^+_l + N^-_l
$$
where
$$
N^+_l = \sum_{t_n \equiv 0 (2)}, \qquad N^-_l = \sum_{t_n \equiv 1 (2)}.
$$
One then finds by the same computation as above for $l=2$ that for $n\geq 2$:
$$
2N_2^-(\spin_{2n+1}) = (2n+1)^g + (2^{2g} -1)(2n+1)^{g-1};
$$
and the left-hand side is replaced by $2N^-_4$ in the case $n=1$.

In \S5 below we shall observe that this is the number of level $2n+1$ theta
functions on the Jacobian and on all the Pryms of a curve of genus $g$.

\space

\beginsection \S4 The height of a representation

\medskip

Suppose, as in the previous section, that the group $G$ of \S1 is simple and
simply-connected.
Given a representation $G\rightarrow SL(V)$ one can associate in a natural way
a level $m_V \in \z$ with which to compute the Verlinde formula of the previous
section. This number---which it seems convenient to refer to as the {\it
height}---has some inte
resting properties and can be defined in several ways (see [KNR]).

Fix a maximal torus in $G$ with Lie algebra $\h$; and let $\Lambda \subset
\h^*$ be the weight lattice, and
$\z[\Lambda]$ the group ring of $\Lambda$. Let
$e:\Lambda\hookrightarrow \z[\Lambda]$ be the natural inclusion. So $e(\lambda
+\mu) = e(\lambda)e(\mu)$. For any $\lambda \in \Lambda$ let $n_{\lambda} \in
\z$ be the multiplicity of $\lambda$ as a weight of the representation $V$.
This determines a natu
ral homomorphism $ch$ of the representation ring of $G$ into $\z[\Lambda]$ by
$$
ch(V) = \sum_{\lambda \in \Lambda} n_{\lambda} e(\lambda) \in \z[\Lambda],
$$
called the {\it formal character} of $V$.

Now let $\cotheta$ be the highest coroot of $G$ with respect to the Cartan
subalgebra $\h$ and normalised Killing form $\< \ ,\ \>$. Then we define:
$$
m_V = {1\over 2} \sum_{\lambda \in \Lambda} n_{\lambda} \< \lambda,\cotheta
\>^2.
\leqno(4.1)
$$
For equivalent definitions see [KNR].
The fundamental result of [F], [KNR] is now the following.

\proclaim (4.2) Theorem. For any simple, simply-connected complex Lie group
$G$, finite dimensional representation $V$, and integer $k$, we have
$$
\dim H^0(\mm(G),\Theta(V)^{\otimes k}) = N_{km_V}(G).
$$

\noindent {\it Remarks.} (i) Of course, although not explicit in the notation,
both sides depend upon the genus $g$.

(ii) Strictly speaking (4.2) is proved only for the classical groups and $G_2$.
See [F].

\medskip

Our present interest in the height invariant is the following example.

\proclaim (4.3) Lemma. If $G =\spin_m$ and $V \cong \c^m$ is the orthogonal
representation, then
$$
m_V =
\cases{
4 & if $m=3$, \cr
2 & if $m\geq 5$. \cr}
$$

\pf In the case $m=3$, $V$ is just the adjoint representation and $m_V = 2h
=4$, where $h$ denotes the dual Coxeter number.
We then have to compute separately the cases $m=6$, $m=2n$ with $n\geq 4$ and
$m=2n+1$ with $n\geq 2$.

First, $\spin_6=SL_4$ and in this case $V \cong \bigwedge^2 \c^4$.
In general, if $G=SL_{n+1}$ then we take the Cartan subalgebra $\h$ to be that
consisting of tracefree diagonal matrices; and let $L_1,\ldots ,L_{n+1}$ be the
standard dual basis in $\h^*$. So the tracefree condition is $L_1+\cdots +
L_{n+1}=0$; and the n
ormalised Killing form is given by:
$$
\langle L_i,L_j \rangle =
\cases{n/(n+1) &if $i=j$,\cr
       -1/(n+1) &if $i\not= j$.\cr
       }
$$
Then $\theta = \cotheta = L_1 -L_{n+1}$; and
one finds that
$V =  \bigwedge^2 \c^{n+1}$ has formal character
$$
ch(V) = \sum_{1\leq i<j\leq n+1} e(L_i +L_j).
$$
Using the above expression for $\<\ ,\ \>$ one finds $m_V = n-1$. Thus for
$\spin_6 $ we obtain $m_V = 2$.

For $\spin_{2n}$, $n\geq 4$, the formal character of $V$ is
$$
ch(V) = \sum_{i=1}^n (e(L_i) + e(-L_i))
$$
where $L_1,\ldots ,L_n$ can be taken to be an {\it orthonormal} basis of
$\h^*$;
whilst $\cotheta = L_1 +L_2$, and so
$$m_V = \sum_{i=1}^n \<L_i,L_1+L_2 \>^2 =2.
$$

Similarly, for $\spin_{2n+1} $ we have
$$
ch(V) = 1 + \sum_{i=1}^n (e(L_i) + e(-L_i)),
$$
for an orthonormal basis $L_1,\ldots ,L_n$,
and $\cotheta =L_1 +L_2$, giving $m_V=2$ as in the even case.

\space

One would like to extend these results for the spin groups to the `twisted'
moduli space $\mm^-(\spin_m)$ defined in \S2, generalising the twisted Verlinde
formula of Thaddeus [T] in the case $m=3$. One obtains a line bundle $\Theta
(V)$ on each of $\mm(\
spin_m)$ and $\mm^-(\spin_m)$ for any representation $V$ of the special
Clifford group.
The appropriate generalisation of the twisted $\spin_3$ Verlinde formula is
then:

\proclaim (4.4) Conjecture. $\dim H^0(\mm^-(\spin_m),\Theta(V)) = -N^+_{m_V}
+N^-_{m_V}$ for $m$ odd, where $N^+$ and $N^-$ are as defined at the end of the
previous section.

We shall elaborate further on this conjecture in \S6 (see also [OW], conjecture
(5.2), to which it is equivalent).
Our present interest is the following identity which follows easily from the
computations of the previous section:
$$
\eqalign{
\dim H^0(\mm^-(\spin_{2n+1}) ,\Theta(\c^{2n+1})) &=
\cases{-N_4^+ +N_4^- & if $n=1$,\cr
   -N_2^+ +N_2^- & if $n \geq 2$;\cr
}\cr
&= 2^{2g-1}(2n+1)^{g-1} + n(2n+1)^{g-1} -2^{2g-1}.\cr}
\leqno(4.5)
$$

\space

\beginsection \S5 Theta functions on Prym varieties

\medskip

We first recall the usual notation for Prym varieties (see [ACGH]).
 For each nonzero half-period $\eta \in J_2(C) \backslash
\{\oo\}$ we have an unramified double cover
$$p : \ctil \rightarrow C.$$
Writing $\jtil = J^{2g-2}(\ctil)$ we have
$$\nmp ^{-1}(K_C) = \peta \cup \peta^- \subset \jtil;$$
where $\peta,\ \peta^-$ are disjoint translates of the same abelian subvariety,
characterised by the condition that for $L\in \nmp^{-1}(K_C)$:
$$
h^0(\ctil, L) \equiv \cases{0& mod 2 if $L\in \peta$,\cr
                              1& mod 2 if $L\in \peta^-$.\cr}
$$
Then the Prym variety of the covering is defined to be $\peta$ (at least in
this section---in the next section it will be more convenient to take the Prym
variety in degree zero).
We shall denote by $\xieta$ the line bundle representing the canonical
principal polarisation on $\peta$,
defined by $2 \xieta = \peta \cap \widetilde \Theta$, where $\widetilde \Theta$
is the theta-divisor in $\jtil$.
Finally, we shall allow also $\eta = 0$ by setting $(P_0 ,\Xi_0) =
(J^{g-1}(C),\Theta)$.

Now let $(A,\Xi)$ be any principally polarised abelian variety of dimension
$g$; and let
$$
H^0(A, m\Xi) = H^0_+(A, m\Xi)\oplus H^0_-(A, m\Xi)
$$
be the decomposition into $\pm$-eigenspaces under the canonical involution of
$A$. Then by writing down a suitable basis of theta functions one can easily
verify that:
$$
\dim H^0_{\pm}(A, m\Xi) =
\cases{
{(m^g \pm 2^g )/ 2} & if $m\equiv 0$ mod 2,\cr
{(m^g \pm 1 )/ 2} & if $m\equiv 1$ mod 2.\cr}
$$

This paper is motivated by the following observation. If we direct sum all {\it
even} theta functions of odd level $m=2n+1$ over all the Pryms $\peta$
(including the Jacobian) we obtain a vector space of dimension
$$
{(2n+1)^g +1 \over 2} +(2^{2g}-1) {(2n+1)^{g-1} +1 \over 2}
$$$$
\eqalign{
&= 2^{2g-1}(2n+1)^{g-1} + n(2n+1)^{g-1} +2^{2g-1}\cr
&=
\cases{
N_2(\spin_{2n+1}) & if $n\geq 2$ \cr
N_4(\spin_3) & if $n=1$.\cr}
\cr}
$$
 From theorem (4.2) and lemma (4.3) it therefore follows that
$$
h^0(\mm(\spin_{2n+1}), \Theta(\c^{2n+1})) = \sum_{\eta \in J_2(C)}
h^0_+(\peta,(2n+1)\xieta),
\leqno(5.1)
$$
for all $n\geq 1$, where $V$ is the orthogonal representation of
$\spin_{2n+1}$.
{\it And similarly,} it follows from (4.5) (assuming the conjecture
(4.4)---which is at least true when $n=1$) that
$$
h^0(\mm^-(\spin_{2n+1}), \Theta(\c^{2n+1})) = \sum_{\eta \in J_2(C)}
h^0_-(\peta,(2n+1)\xieta).
\leqno(5.2)
$$

\noindent
{\bf (5.3) Remark.} It may be noted from the remarks following (3.2) (and again
using (4.3)) that among the even spin groups $\spin_8$ alone exhibits an
analogous numerology, namely:
$$
h^0(\mm(\spin_{8}), \Theta(\c^{8})) = \sum_{\eta \in J_2(C)}
h^0(\peta,8\xieta).
$$

\medskip

The identities (5.1) and (5.2) were explained for $n=1$ in [O]; the purpose of
the remainder of this section is to account similarly for $n > 1$; that is, to
prove the generalisation (0.3) of (0.1) stated in the introduction.
In so doing we shall give only briefly those parts of the construction for
which the general case differs only trivially from the case $n=1$. At such
points we shall refer to [O] for the details.

\space

It will be convenient to denote the two-component variety $\mm(\spin_n) \cup
\mm^-(\spin_n)$ by $\nn(n)$; and the theta line bundle $\Theta(\c^n)$ by
$\Theta_{\nn}$. We shall construct a homomorphism
$$
H^0(\nn(n),\Theta_{\nn})\dual \rightarrow \bigoplus_{\eta \in J_2(C)}
H^0(P_{\eta},n \xieta),
$$
one summand at a time; and then check that the parity corresponds on each side
as in (0.3).
By the K\"unneth formula this is equivalent---up to the assertion about
parity---to finding, for each $\eta \in J_2(C)$, a naturally occuring divisor
$$
\dd_{\eta} \subset P_{\eta} \times \nn(n)
$$
which is in the linear system $|n\xieta + \Theta_{\nn}|$. (Here, of course,
$\xieta$ and $\Theta_{\nn}$ denote by abuse of notation the respective line
bundles pulled back to the product.)

\medskip
\noindent
{\bf The case $\bf \eta =0$.}
\medskip

In this case, which we shall consider first, $P_{0}=J^{g-1}(C)$ and $\Xi_0
=\Theta_J$ is the usual theta divisor.

We consider the morphism of varieties
$$
\alpha : \nn(n) \rightarrow \mm(SC_n) \rightarrow \mm(SL_n)
$$
induced by the group homomorphisms
$$
SC_n \rightarrow SO_n \hookrightarrow SL_n.
$$
In other words, to each Clifford bundle $E$ we associate the vector bundle
$E(\c^n)$; then, as already mentioned in \S1, semistability of $E$ implies that
of $E(\c^n)$ ([R2], proposition 2.17).
This assignment is functorial with respect to base change for families, so by
the coarse moduli property we obtain the above morphisms.

Now we define $\dd_0$ to be the pull-back via
$$
J^{g-1}(C) \times \nn(n) \map{id \times \alpha} J^{g-1}(C) \times \mm(SL_n)
\map{\otimes} \mm(n,n(g-1))
$$
(where, of course, we are identifying $\mm(SL_n)$ with the moduli space of
semistable vector bundles of rank $n$ and trivial determinant, and the second
map is tensor product of vector bundles)
of the canonical theta divisor (see \S1)
$$
\Theta_{n,n(g-1)} = \{V| H^0( V) \not= 0\} \subset \mm(n,n(g-1)).
$$

It follows from the discussion of \S1 that $\dd_0$ belongs to the linear system
$|n\Theta_J + \Theta_{\nn}|$: first of all, the pull-back of
$\Theta_{n,n(g-1)}$ to $J^{g-1}(C) \times \mm(SL_n)$ restricts on a fibre
$\{L\} \times \mm(SL_n)$ as the restrict
ion of $\Theta_{n,0}$ from $\mm(n,0)$---this is by definition of
$\Theta_{n,0}$---and we have seen that this is just $\Theta(\c^n)$. Hence the
pull-back to $J^{g-1}(C) \times \nn(n)$
is $\Theta_{\nn}$ on fibres $\{L\} \times \nn(n)$.

On the other hand, restriction to fibres $J^{g-1}(C) \times \{V\}$, for any
vector bundle $V$ with rank $n$ and trivial determinant, is independent of $V$
by $(1.5')$. We may therefore take $V$ to be the trivial bundle of rank $n$.
Then it is clear from the usual determinantal description of $\Theta_J$ that
$\Theta_{n,n(g-1)} $ restricts to $n\Theta_J$.

\medskip
\noindent
{\bf The case $\bf \eta \not= 0$.}
\medskip

To begin, it is necessary to note that semistability of a bundle is preserved
under pull-back.

\proclaim (5.4) Lemma.
Let $p : \cctil \rightarrow C$ be any unramified cover of smooth projective
curves.  Then, if a vector bundle $V\rightarrow C$ is semistable then
$p^*V\rightarrow \cctil$ is semistable.

\pf By the Narasimhan-Seshadri theorem $V$ is induced from a projective unitary
representation of the fundamental group $\pi_1(C)$. Since $\cctil$ is an
unramified cover its fundamental group injects into $\pi_1(C)$, and the
restriction of the above repre
sentation then induces the pull-back bundle, which is consequently semistable.

\proclaim (5.5) Corollary. Let $p : \cctil \rightarrow C$ be as in lemma (5.4),
and $E\rightarrow C$ a semistable $G$-bundle. Then $p^*E \rightarrow \cctil$ is
semistable.

\pf The same argument as in the above proof works for $G$-bundles by
Ramanathan's generalisation of the Narasimhan-Seshadri theorem [R1];
alternatively apply (5.4) to the adjoint bundle $\ad E$: by [R2], corollary
2.18, semistability of $E$ is equivalent
to semistability of $\ad E$ as a vector bundle.

\medskip

Let us now return to the double cover $p:\ctil \rightarrow C$. Noting that for
a Clifford bundle $E \rightarrow C$ the spinor norm satisfies $\nm ( p^*E) =
p^* \nm ( E)$, and this has even degree, it follows from (5.5) that we obtain a
morphism of moduli
spaces
$$
u = p^*: \nn_C(n)  \rightarrow \mm_{\ctil}(\spin_n).
$$

\proclaim (5.6) Proposition. $u^*\Theta_{\mm_{\ctil}}(\c^n) = 2\Theta_{\nn}$.

{\it Proof.}
Let $E\rightarrow C\times S$ be an arbitrary family of semistable
$SC_n$-bundles, and let $F=E(\c^n)$ be the associated family of vector bundles
via the orthogonal representation. Let $\ftil = (p\times {\rm id})^*F$ be the
pull-back of the family by the d
ouble cover:
$$
\matrix{
\ctil \times S&&\map{p\times {\rm id}}&&C\times S\cr
&&&&\cr
&\hidewidth \pitil \searrow &&\swarrow \pi \hidewidth&\cr
&&&&\cr
&&S.&&\cr}
$$

It is clear from the discussion of \S1 that to prove the proposition it
suffices to show that
$$
\Theta(\ftil) = 2 \Theta (F):
$$
i.e. the line bundle $\Theta(\c^n) \rightarrow \mm_C(SC_n)$ represents the
functor $E\mapsto \Theta (F)$, while the line bundle
$u^*\Theta_{\mm_{\ctil}}(\c^n) \rightarrow \mm_C(SC_n)$ represents the functor
$E\mapsto \Theta(\ftil)$.

So to compute $\Theta(\ftil)$, first note that by the projection formula
applied to $p\times {\rm id}$ we have, for $i=0,1$:
$$
\eqalign{
R^i_{\pitil}(\ftil) &= R^i_{\pi}(F\otimes p_* \oo_{\ctil}) \cr
                    &= R^i_{\pi}(F\oplus F\otimes \eta)\cr
                    &=R^i_{\pi}(F) \oplus R^i_{\pi}(F\otimes \eta).\cr}
$$
If we fit the direct images $R^i_{\pi}(F)$ into an exact sequence (1.3), and
$R^i_{\pi}(F\otimes \eta)$ into a similar sequence with middle terms ${K^0}'
\map{\phi'} {K^1}'$, then we get an exact sequence:
$$
0\rightarrow R^0_{\pitil}\ftil
\rightarrow K^0 \oplus {K^0}'
\map{
{\rm diag}(\phi,\phi')
} K^1 \oplus {K^1}'
\rightarrow R^1_{\pitil}\ftil\rightarrow 0.
$$
It follows at once that
$$
{\rm Det}(\ftil) = {\rm Det}(F)\otimes {\rm Det}(F\otimes \eta).
$$
But since the bundle $F$ has trivial determinant we can replace Det by $\Theta$
here. And since $\Theta(F\otimes \eta) = \Theta(F)$ by corollary (1.6), we
obtain $\Theta(\ftil) = 2 \Theta(F)$ as required.

\medskip

As a consequence of proposition (5.6)     , we see that the pull-back via the
map
$$
 {\rm incl} \times u : \peta \times \nn(n) \rightarrow \jtil \times
\mm_{\ctil}(\spin_n),
$$
of the divisor $\widetilde \dd_0$ constructed (with $C$ replaced by $\ctil$)
for the case $\eta =0$, is a divisor
$$
\ee_{\eta} \in |2n \xieta + 2\Theta_{\nn}|.
$$

\proclaim (5.7)    Proposition. $\ee_{\eta}$ has multiplicity 2, i.e.
$\ee_{\eta} = 2 \dd_{\eta}$ for a divisor $\dd_{\eta} \in |n \xieta +
\Theta_{\nn}|$.

For the proof of this we refer to [O], proposition (2.6), where the case $n=3$
is proved. The proof for the general case is not essentially different except
in one respect: for $n=3$ there exists a {\it universal} orthogonal bundle
$E(\c^3) = \ad E \rightarrow C\times \nn(3)$, whereas we cannot expect this in
the general case.
However, since the result to be proved is local, it suffices to prove the
result on spaces $\peta \times S$ where $S$ is the base of an arbitrary family
of semistable Clifford bundl
es, and then apply the universal moduli property of $\nn(n)$. The crucial point
in both cases is that the bundles $E(\c^n)$ (or $\ad E$ when $n=3$) are {\it
self-dual}: see [O] lemma (2.11).

\medskip

One can view the divisor $\dd_{\eta}$ as giving a rational map
$$
\eqalign{
g_{\eta }: \nn(n) & \rightarrow |n \xieta | \cr
                E  &\mapsto  \{ L \in \peta | H^0(\ctil , L\otimes p^*E(\c^n))
\not= 0 \}.\cr}
$$
It is easy to see by Riemann-Roch, Serre duality and the fact that $E(\c^n)$ is
self-dual, that $g_{\eta}(E)$ is a {\it symmetric} divisor (see [O], \S3).
This means that each component of $\nn(n)$ maps either into $|n\xieta |_+ =
\p H^0_+(\peta , n \xieta)$ or into  $|n\xieta |_- = \p H^0_-(\peta , n
\xieta)$; and the claim is:

\proclaim (5.8)   Proposition. For each $\eta \in J_2(C)$, $g_{\eta}$ maps
$\mm(\spin_n)$ into $|n\xieta |_+$ and $\mm^-(\spin_n)$ into $|n\xieta |_-$.

As a consequence we obtain the homomorphisms of (0.3) in the introduction,
$$
f_{\eta}^{\pm} =
(g_{\eta}^*)\dual : H^0(\mm^{\pm}(\spin_n), \Theta(\c^n))\dual
 \rightarrow
H^0_{\pm}(\peta, n\xieta)
\leqno (5.9)
$$
dual to the pull-back of hyperplane sections.

The proof of (5.8)   is identical to that of [O], \S3, if we replace the
congruence [O],(3.7) with its generalisation in
proposition (2.11) of the present paper. This implies that for $E\in
\mm(\spin_n)$ the divisor $g_{\eta}(E) \subset \peta $ passes through all odd
theta-characteristics, while for $E\in \mm^-(\spin_n)$ it passes through all
even theta-characteristics. By
[O], lemma (3.6) the odd theta-characteristics are the base points of the
linear subsystem $|n\xieta |_+$ and the even theta-characteristics are the base
points of the linear subsystem $|n\xieta |_- $; so the result follows.

\space

\beginsection \S6 Spin reciprocity

\medskip

We are not going to discuss, in this article, the question of whether the
homomorphisms (5.9) are isomorphisms as one should expect. However, in this
final section we wish to fit this question into a broader context by giving an
alternative point of view
on the relationship between Prym varieties and spin bundles. We consider
complex $O_2$- and $\pin_2$-bundles on the curve $C$; and we shall show that
these can be identified with anti-invariant line bundles over double covers of
$C$ ((6.6) and (6.7)).

First note that $O_2$ and $\pin_2$ are isomorphic groups; each is generated by
$SO_2 \cong \spin_2 \cong \c^*$ and a fixed element $\sigma$ of order 2:
$$
O_2 \cong \pin_2 \cong \<\c^*, \sigma | \sigma^2 = 1,\ \sigma z = z^{-1}\sigma
\ {\rm for}\ z\in \c^*\>.
\leqno(6.1)
$$
The double cover of $O_2$ by $\pin_2$ is then induced by the squaring map on
$\c^*$:
$$
\matrix{
1&\rightarrow&\c^{*}&\rightarrow&\pin_2&\map{\nu}&\z/2&\rightarrow&0\cr
&&&&&&&&\cr
&&\hidewidth sq\downarrow\ \ &&\downarrow&&\Vert&&\cr
&&&&&&&&\cr
1&\rightarrow&\c^*&\rightarrow&O_2&\map{\nu}&\z/2&\rightarrow&0.\cr}
\leqno(6.2)
$$
We shall use the same letter $\nu$ for the map induced in cohomology
$$
\nu : H^1(C,\oo O_2) \rightarrow H^1(C,\z/2) \cong J_2(C).
\leqno(6.3)
$$
The first space $H^1(C,\oo O_2)$ parametrises isomorphism classes of
$O_2$-bundles on $C$; let us represent such a bundle by a Cech cocycle $g =
\{g_{ij}\}$ with respect to some open cover $\{U_{ij}\}$ of $C$. Each $g_{ij} =
\eta_{ij} \gtil_{ij}$ where $\
eta_{ij} \in \{1,\sigma \}$ and $\gtil_{ij} \in \oo^*(U_i \cap U_j)$. It will
be convenient to identify the groups $\{1,\sigma \} \cong \{\pm 1 \} \cong
\z/2$. Then we observe, first, that $\{\eta_{ij} \}$ is a cocycle representing
the image $\eta = \nu (
g)$ under (6.3).

Second, it is easy to verify that the  cocycle conditions satisfied by
$\{g_{ij}\}$, together with the relations (6.1), imply:
$$
\gtil_{ij}^{\eta_{ji}}\gtil_{ji} = 1, \qquad
\gtil_{ij}^{\eta_{jk}\eta_{ki}}\gtil_{jk}^{\eta_{ki}}\gtil_{ki} = 1.
\leqno(6.4)
$$
Although $\{\gtil_{ij}\}$ does not, therefore, define a cocycle on $C$, we
claim that it {\it does} define one on the unramified double cover $\ctil
\rightarrow C$ associated to $\eta = \{\eta_{ij} \}$.

To see this, we first construct an open cover $\{U_i^+,U_i^- \}$
 of $\ctil$, where $U_i^+ \cong U_i^- \cong U_i \subset C$, by glueing
$$
U_i^{\pm} \rightarrow  U_j^{\pm \eta_{ij}}
$$
over the intersection $U_i \cap U_j$. (This is one way to construct the curve
$\ctil$.) We then have a collection of functions
$$
\{\gtil_{ij} \in \oo^*(U_i^+ \cap U_j^{\eta_{ij}}),\quad
\gtil_{ij}^{-1} \in \oo^*(U_i^- \cap U_j^{-\eta_{ij}})
\}.
\leqno(6.5)
$$
Notice that switching $i$ and $j$ replaces $\gtil_{ij}$ by
$\gtil_{ji}^{\eta_{ij}}$; thus one sees that the relations (6.4) are precisely
the cocycle conditions for the collection (6.5); {\it i.e. the functions (6.5)
are transition functions for an anti-i
nvariant line bundle on the curve $\ctil$.}

If $\eta = 0 \in J_2(C)$ then the double cover is trivial, i.e. two copies of
$C$, and an anti-invariant line bundle just means an arbitrary line bundle on
$C$ (together with its inverse on the other component). If $\eta \not= 0$ then
anti-invariant line
bundles are parametrised by
$$
\ker \nm_{p_{\eta}} \subset J(\ctil).
$$
(Note that such a line bundle necessarily has degree zero.) In summary,
therefore, we have shown:

\proclaim (6.6) Proposition. $$H^1(C,\oo O_2) \cong \pic(C) \cup \bigcup_{\eta
\in J_2(C)\backslash \{0\}}\ker \nm_{p_{\eta}}.$$

Recall that for each $\eta \in J_2(C)\backslash \{0\}$ the kernel of the norm
map has two connected components:
$$
\ker \nm_{p_{\eta}} = \prym(C,\eta) \cup \prym^-(C,\eta)
$$
where $\prym(C,\eta) = \im(1 -\sigma)$ and where $\sigma$ denotes the
endomorphism of $J(\ctil)$ induced by the sheet-interchange over $C$.
(Compare \S5: $\prym(C,\eta)$ is, of course, a translate of $\peta$.)

\proclaim (6.7) Proposition. The subvariety of $H^1(C,\oo O_2)$ parametrising
bundles which lift to a $\pin_2$-bundle on $C$ is precisely
$$
\pic^{\rm even}(C) \cup \bigcup_{\eta \in J_2(C)\backslash \{0\}}\prym(C,\eta).
$$

\pf An $O_2$-bundle lifts to a $\pin_2$-bundle if and only if its cohomology
class is in the image of the squaring map (6.2); which is equivalent to saying
that the line bundle we obtain on the double cover has a square root.
So if $\eta =0$ then the line bundle we obtain on $C$ has even degree; whilst
if $\eta \not= 0$ then we obtain a line bundle $L = N^2 \rightarrow \ctil$,
where $\sigma(N) = N^{-1}$. But this implies that $L= N\otimes \sigma (N)^{-1}$
i.e. that $L \in \im(
1 - \sigma) = \prym(C,\eta)$.

\space

`Spin reciprocity' is the term used in [OW] to refer to the fact that for odd
numbers $l,m \geq 5$,
$$
N_l^-(\spin_m) = N_m^-(\spin_l), \qquad N_{2l}^-(\spin_3) = N_3^-(\spin_l).
\leqno(6.8)
$$
This can be interpreted by considering the variety (as in \S5)
$$
\nn(m) = \cases{ \mm(\spin_m) \cup \mm^-(\spin_m) & for odd $m\geq 3$,\cr
                 J_2(C) & for $m=1$.\cr}
$$
(The convention for $\nn(1)$ is included mainly for numerological convenience;
though one may note that in each case $\nn(m)$ is characterised as being the
canonical Galois cover of $\mm(SO_m)$ with fibre $J_2(C)$.)

Now suppose that there exists on $\nn(m)$ a `Pfaffian' line bundle $\ll$ such
that $\ll^2 = \Theta(\c^m)$. From the construction of [KNR] it would follow
that, more generally, for any representation $V$ of $SC_m$ we have $\Theta(V) =
\ll^{m_V}$.
In fact the work of Sorger [So] suggests that $\ll$ should exist in general
only as a Weil divisor class, Cartier over the stable points of the moduli
space; but in this case one might still expect that the statements below about
$H^0$s are unaffected.

Now conjecture (4.4) implies that (for $m$ odd, of course)
$$
\dim H^0(\nn(m), \ll^l) = 2N^-_l(\spin_m),
$$
and so (6.8) becomes
$$
\dim H^0(\nn(m), \ll^l) = \dim H^0(\nn(l),\ll^m)
\qquad
\hbox{for $l,m$ both odd.}
\leqno(6.9)
$$

One may now ask, on the one hand, whether there is a natural duality between
these two vector spaces, and on the other hand how the {\it even} spin groups
fit into this picture. For the latter question, in particular, the only
numerological clue appears t
o be remark (5.3); though one can expect $\spin_8$ to be special on account of
triality. However,
the point that we wish to make here is that if we take
$$
\nn(2) = J(C) \cup \bigcup_{\eta \in J_2(C)\backslash \{0\}} \prym (C,\eta)
$$
with line bundle $\ll$ restricting to the principal polarisation on each
component, then by (5.1) and (5.2) the reciprocity relation (6.9) remains true
if one of $l,m$ equals 2---and indeed a pairing is constructed via
(5.9)---whilst by (6.7) the variety
$\nn(2)$ has a natural interpretation as a moduli space of $\pin_2$-bundles on
$C$.

\space

\vfill\eject

{\capit REFERENCES}

\bigskip
\item{[ACGH]}  {\capit E. Arbarello, M. Cornalba, P.A. Griffiths, J. Harris}:
Geometry
of Algebraic Curves, Springer, 1985;
\medskip
\item{[B1]} {\capit A. Beauville}: Fibr\'es de rang 2 sur une courbe, fibr\'e
d\'eterminant et fonctions theta, Bull. Soc. Math. France, 116 (1988) 431--448;
\medskip
\item{[B2]} \quad ---\quad : --- II, Bull. Soc. Math. France, 119 (1991)
259--291;
\medskip
\item{[B3]} \quad --- \quad : Conformal blocks, fusion rules and the Verlinde
formula, preprint 1994;
\medskip
\item{[Bo]} {\capit N. Bourbaki}: Alg\`ebre II chap. 10, Formes
sesquilin\'eaires et formes quadratiques, Hermann 1959;
\medskip
\item{[DN]} {\capit J.-M. Drezet, M.S. Narasimhan}: Groupe de Picard de
vari\'et\'es
de modules de fibr\'es semi-stable sur les courbe alg\'ebriques, Invent. Math.
97 (1989) 53--94;
\medskip
\item{[F]} {\capit G. Faltings}: A proof for the Verlinde formula, J. Alg.
Geometry 3 (1994) 347--374;
\medskip
\item{[K]} {\capit V. Kanev}: Spectral curves, simple Lie algebras and
Prym-Tyurin varieties, Proc. Symp. Pure Math. 49 (1989) Part 1, 627--645;
\medskip
\item{[KNR]} {\capit S. Kumar, M.S. Narasimhan, A. Ramanathan}: Infinite
Grassmannian and moduli space of $G$-bundles, preprint 1993;
\medskip
\item{[O]} {\capit W.M. Oxbury}: Anticanonical Verlinde spaces as theta spaces
on Pryms, preprint 1993;
\medskip
\item{[OW]} {\capit \quad --- \quad, S.M.J. Wilson}: Reciprocity laws in the
Verlinde formulae for the classical groups, preprint 1994;
\medskip
\item{[R]} {\capit S. Ramanan}: Orthogonal and spin bundles over hyperelliptic
curves, in \lq Geometry and Analysis', papers dedicated to V.K. Patodi,
Springer (1981);
\medskip
\item{[R1]} {\capit A. Ramanathan}: Stable principal bundles on a compact
Riemann surface, Math. Ann. 213 (1975) 129--152;
\medskip
\item{[R2]} {\capit \quad --- \quad}: Stable principal bundles on a compact
Riemann surface: construction of moduli space, Bombay thesis, 1976;
\medskip
\item{[S]} {\capit J.-P. Serre}: Rev\^etements \`a ramification impaire et
th\^eta-caract\'eristiques, C. R. Acad. Sci. Paris 311(I) (1990) 547--552;
\medskip
\item{[So]} {\capit C. Sorger}: Groupe de Picard de la vari\'et\'e des
th\^eta-caract\'eristiques des courbes planes, to appear in the proceedings of
the 1993 Europroj conference at Catania;
\medskip
\item{[T]} {\capit M. Thaddeus}: Conformal field theory and the cohomology of
the moduli space of stable bundles, J. Diff. Geometry 35 (1992) 131--149.

\space

{\addressit
\ni Department of Mathematical Sciences,

\ni Science Laboratories,

\ni South Road,

\ni Durham DH1 3LE,

\ni U.K.}

\medskip

\ni {\addressit E-mail:} {\eightrm w.m.oxbury@durham.ac.uk }

\end